\documentstyle[11pt]{article}

\makeatletter

\@addtoreset{equation}{section}
\makeatother
\def\C{\mbox{\boldmath $C$}}

\def\<{\langle}
\def\>{\rangle}

\newtheorem{lem}{Lemma}[section]
\newtheorem{theo}{Theorem}[section]
\newtheorem{rem}{Remark}[section]

\makeatletter
   
   \@addtoreset{equation}{section}
\makeatother

\begin{document}
\title{\bf Revisit on How to Derive Asymptotic Profiles to Some Evolution Equations}
\author{Ryo IKEHATA\thanks{Corresponding author: ikehatar@hiroshima-u.ac.jp} \\ {\small Department of Mathematics, Graduate School of Education, Hiroshima University} \\ {\small Higashi-Hiroshima 739-8524, Japan}}
\maketitle
\begin{abstract}
We consider the Cauchy problem in ${\bf R}^{n}$ for heat and damped wave equations. We derive asymptotic profiles to those solutions with weighted $L^{1,1}({\bf R}^{n})$ data by presenting a simple method. 
\end{abstract}

\section{Introduction}
\footnote[0]{Keywords and Phrases: Heat equation; Damped wave equation; Fourier analysis; Asymptotic profiles; Low frequency; High frequency; Weighted $L^{1}$-initial data.}
\footnote[0]{2000 Mathematics Subject Classification. Primary 35L60, 37L15; Secondary 35L70, 35B40.}
We first consider the Cauchy problem for heat equations in ${\bf R}^{n}$ ($n \geq 1$):
\begin{equation}
v_{t}(t,x) - \Delta v(t,x) = 0,\ \ \ (t,x)\in (0,\infty)\times {\bf R}^{n} ,\label{eqn}
\end{equation}
\begin{equation}
v(0,x)= v_{0}(x), \quad x\in {\bf R}^{n} ,\label{initial}
\end{equation}
where the initial datum $u_{0}$ is taken from the weighted $L^{1}$-space:
\[v_{0} \in L^{2}({\bf R}^{n}) \cap L^{1,1}({\bf R}^{n}),\]
where
\[f \in L^{1,k}({\bf R}^{n}) \Leftrightarrow f \in L^{1}({\bf R}^{n}),\,\Vert f\Vert_{1,k} := \int_{{\bf R}^{n}}(1 + \vert x\vert)^{k}\vert f(x)\vert dx < +\infty, \quad k \in {\bf N}\cup \{0\}.\]
\noindent
It is known that for each $v_{0} \in L^{2}({\bf R}^{n})$ the problem (1.1)-(1.2) admits a unique solution $v \in C([0,+\infty);L^{2}({\bf R}^{n}) \cap C^{1}((0,+\infty);H^{2}({\bf R}^{n}))$ (cf.\,\cite{p}) \\


Our first target is to give a simple alternative proof of the following known result (cf. \cite{g}).

\begin{theo} Let $n \geq 1$. If $v_{0} \in L^{2}({\bf R}^{n})\cap L^{1,1}({\bf R}^{n})$, then the solution $v(t,x)$ to problem {\rm (1.1)}-{\rm (1.2)} satisfies
\[\Vert v(t,\cdot)-P_{0}G(t,\cdot)\Vert \leq Ct^{-\frac{n}{4}-\frac{1}{2}}(\int_{{\bf R}^{n}}\vert x\vert\vert v_{0}(x)\vert dx),\quad t > 0,\]
where
\[P_{0} := \int_{{\bf R}^{n}} v_{0}(x)dx,\]
\[G(t,x) := \frac{1}{(\sqrt{4\pi t})^{n}}e^{-\frac{\vert x\vert^{2}}{4t}},\]
and $C > 0$ is a constant.
\end{theo}
\begin{rem}{\rm It follows from the recent result due to Said-Houari \cite[Theorem 3.2]{said} which derived a more general result from the viewpoint of $W^{k,p}$-decay estimate of solutions (see also \cite{g}) that
\[\Vert v(t,\cdot)\Vert \leq Ct^{-\frac{n}{4}-\frac{1}{2}}(\int_{{\bf R}^{n}}\vert x\vert\vert v_{0}(x)\vert dx) + C\vert P_{0}\vert t^{-\frac{n}{4}}.\]
This implies that under the condition $P_{0} \ne 0$ one has
\[\Vert v(t,\cdot)\Vert = O(t^{-\frac{n}{4}}), \quad t \to +\infty,\]
even if the initial datum belongs to the weighted $L^{1}$ space. Furthermore, it is well-known that  
\[\Vert G(t,\cdot)\Vert = O(t^{-\frac{n}{4}}), \quad t \to +\infty.\]
These observation and Theorem 1.1 imply that the asymptotic profile of solutions to the heat equation (1.1) becomes so called a multiple of the Gauss kernel as $t \to +\infty$. Of course, this is a well-known fact. Indeed, Karch \cite[Lemma 3.2]{K} has already announced the similar fact that
\[\Vert v(t,\cdot)-P_{0}G(t,\cdot)\Vert = o(t^{-n/4}),\quad t \to +\infty,\]
in the case when $v_{0}$ belongs to the usual $L^{1}({\bf R}^{n})$ space (not weighted $L^{1}$ space!). By restricting the initial datum $v_{0}$ to the weighted $L^{1}$-space, we can find a more precise decay order (see also \cite{g} for the same result as in the case of $v_{0} \in L^{1,1}({\bf R}^{n})$).  Our new point of view is to introduce a new simple proof  of this fact, and our method has a possibility widely applied to some other evolution equations including damped wave equations below.

A similar result to the wave equations with structural damping:
\[u_{tt}-\Delta u + (-\Delta)^{\theta}u_{t} = 0, \quad \theta \in [0,1],\]
will be announced in our forthcoming project. }
\end{rem}

Next let us consider the Cauchy problem for damped wave equations in ${\bf R}^{n}$ ($n \geq 1$):
\begin{equation}
u_{tt}(t,x) - \Delta u(t,x)  + u_{t}(t,x) = 0,\ \ \ (t,x)\in (0,\infty)\times {\bf R}^{n} ,\label{eqn}
\end{equation}
\begin{equation}
u(0,x)= u_{0}(x), \quad u_{t}(0,x) = u_{1}(x), \quad x\in {\bf R}^{n} ,\label{initial}
\end{equation}
where the initial data $u_{0}$ and $u_{1}$ are also taken from the weighted $L^{1}$-space:
\[[u_{0},u_{1}] \in (H^{1}({\bf R}^{n}) \cap L^{1,1}({\bf R}^{n})) \cap (L^{2}({\bf R}^{n}) \cap L^{1,1}({\bf R}^{n})).\]
\noindent
Then we can find that the problem (1.3)-(1.4) admits a unique weak solution\\
 $u \in C([0,+\infty);H^{1}({\bf R}^{n})) \cap C^{1}([0,+\infty);L^{2}({\bf R}^{n}))$.\\

Since Nishihara \cite{N} studied the asymptotic profiles to the quasilinear damped wave  equations, which have its origin in the research of the system of hyperbolic conservation, many mathematicians are producing many interesting papers about the diffusion phenomenon of the damped wave equations.

The $L^{p}$-estimates for the difference $u(t,\cdot)-v(t,\cdot)$ based on the Fourier analysis can be found in \cite{K} and \cite{Na}, where $u$ is the solution to (1.3)-(1.4), and $v(t,x)$ is the solution to (1.1)-(1.2) with $v_{0}(x) = u_{0}(x) + u_{1}(x)$. They solved concretely the ODE after the Fourier transformed equation with parameter $\xi$:
\[\hat{u}_{tt}(t,\xi)+\vert\xi\vert^{2}\hat{u}(t,\xi) + \hat{u}_{t}(t,\xi) = 0\]
to proceed the estimates of transformed solution $\hat{u}(t,\xi)$ in the low and high frequency parameter $\xi$ in the Fourier space. In particular, Karch \cite{K} introduced the following equality (in fact, he treated a more general nonlinear equations):
\[\lim_{t \to +\infty} t^{\frac{N}{2}(1-\frac{1}{p})}\Vert u(t,\cdot)-C_{*}G(1+t,\cdot)\Vert_{p} = 0,\]
where $p \in [1,\infty]$, and 
\[C_{*} := \int_{{\bf R}^{n}}(u_{0}(x) + u_{1}(x)))dx.\]
\noindent
The sharp results about the asymptotic expansion of the solution $u(t,x)$ as $t \to +\infty$ was found in \cite{N-2} in the case when $n = 3$:
\[u(t,\cdot) \approx v(t,\cdot) + e^{-t/2}w(t,\cdot), \quad t \to +\infty,\]
where $w(t,x)$ is a solution to the free wave equation
\begin{equation}
w_{tt}(t,x) - \Delta w(t,x) = 0,
\end{equation}
\[w(0,x)= u_{0}(x), \quad w_{t}(0,x) = u_{1}(x).\]
He (\cite{N-2}) used the fundamental solutions to (1.3) based on the famous text book due to Courant-Hilbert, so the restriction on the dimension $n$ seems to be necessary.
These types of asymptotic expansions of global solutions to (1.3)-(1.4) with a power type of nonlinearity were also developed in \cite{HO} ($n = 2$) , \cite{MN} ($n = 1$) and \cite{HKN}(from the viewpoint of weighted $L^{2}$ data).

On the other hand, the abstract theory in Hilbert spaces about the diffusion phenomenon can be found in \cite{CH}, \cite{IN} and \cite{RTY}.

Recently, Said-Houari \cite{said} derived the shaper decay estimates for the difference on $u-v$ in terms of weighted $L^{1}$ initial data. Under the assumption that $\displaystyle{\int_{{\bf R}^{n}}}u_{j}(x)dx = 0$ ($j = 0,1$) he derived the sharp $L^{\infty}$ estimate of the difference $u-v$ based on the two methods from \cite{HM} and \cite{Ik-3}.

Quite recently, Kawakami-Ueda \cite{KU} reconsidered the "nonlinear" version of the problem (1.3)-(1.4) with the nonlinear term $F(t,x,u)$ in the case when $n = 1,2,3$. Their research is also based on a viewpoint of the weighted $L^{1}$-data, i.e., $u_{0} \in W_{k}^{1,1} \cap L^{1,\infty}$ and $u_{1} \in L_{k}^{1}\cap L^{\infty}$. The restriction on the dimension $n$ comes from their method, in fact, they used a similar framework to the Nishihara's one \cite{N-2} based on several estimates for the fundamental solutions of (1.3).  While, we have to mention to the work due to Hosono \cite{H}, in which the asymptotic behavior of solutions of nonlinear problem for (1.3) was studied around 2006 by the Fourier analysis similar to the method introduced in this paper. However, it should be noted that the method presented in this paper basically constructed in 2003.

The purpose in this paper  is to find the asymptotic profile as $t \to +\infty$ of the solution $u$ to problem (1.3)-(1.4) in terms of the "weighted $L^{1}$-initial data" based on an idea to derive Theorem 1.1 above. That idea has its origin in \cite[Lemma 3.1]{Ik-3} and \cite[Lemma 2.3]{Ik-2}. The viewpoint from the weighted initial data seems not so new as is already mentioned (see \cite{Ik-3}, \cite{KU} and \cite{said}). Our novelty is to introduce a simple new method in the case when we derive asymptotic profiles, and our argument is independent from the restriction on the dimension $n$. The term "simple" means that we have only to observe the Fourier transformed initial data thoroughly in order to catch the asymptotic profiles of solutions, that is, the asymptotic state is determined by a decomposition of the Fourier transformed "initial data" (see (3.8)). By this idea we can also consider the higher order expansions of solutions together with applications to the other type of evolution equations, but these applications will be announced in a series of forthcoming projects. Unfortunately, at present our method seems not to be applied to the nonlinear case as in \cite{KU}.\\

Our main target is to give a simple proof of the following fact from the view point of the weighted $L^{1,1}$-initial data.

\begin{theo} Let $n \geq 1$. If $[u_{0},u_{1}] \in (H^{1}({\bf R}^{n}) \cap L^{1,1}({\bf R}^{n})) \cap (L^{2}({\bf R}^{n}) \cap L^{1,1}({\bf R}^{n}))$, then the solution $u(t,x)$ to problem {\rm (1.3)}-{\rm (1.4)} satisfies
\[\Vert u(t,\cdot)-(P_{00}+P_{01})G(t,\cdot)\Vert \leq Ct^{-\frac{n}{4}-\frac{1}{2}}(\Vert u_{0}\Vert_{1,1}+\Vert u_{1}\Vert_{1,1}+\Vert u_{0}\Vert + \Vert\nabla u_{0}\Vert + \Vert u_{1}\Vert),\quad t > 0,\]
where
\[P_{00} := \int_{{\bf R}^{n}} u_{0}(x)dx,\quad P_{01} := \int_{{\bf R}^{n}} u_{1}(x)dx.\]
\end{theo}
In 2003 Ikehata \cite{Ik-3} proved the following result for the solution $u(t,x)$ to problem (1.3)-(1.4) based on the previously computed one due to Matsumura \cite{Ma}:
\begin{equation}
\Vert u(t,\cdot)\Vert \leq C I_{0}(1+t)^{-\frac{n}{4}-\frac{1}{2}} + C\vert P_{00} + P_{01}\vert(1+t)^{-\frac{n}{4}},
\end{equation}
where
\[I_{0} := \Vert u_{0}\Vert_{H^{1}}+\Vert u_{0}\Vert_{1,1} + \Vert u_{1}\Vert +\Vert u_{1}\Vert_{1,1}.\]
This implies that in the case when $\vert P_{00} + P_{01}\vert \ne 0$, we have at most
\[\Vert u(t,\cdot)\Vert = O(t^{-\frac{n}{4}}), \quad t \to +\infty.\]
Furthermore, it follows from the same observation as in Remark 1.1 that 
\[\Vert G(t,\cdot)\Vert = O(t^{-\frac{n}{4}}), \quad t \to +\infty.\]
So, the result in Theorem 1.2 implies that in the case when $P_{00} + P_{01} \ne 0$ the asymptotic profile of the solution $u(t,x)$ to problem (1.3)-(1.4) as $t \to +\infty$ becomes a multiple of the Gauss kernel.  The result in Theorem 1.2 becomes an improvement from the viewpoint of the $L^{1,1}$-initial data (cf. \cite{CH}, \cite{IN}, \cite{K}, \cite{Na}, \cite{N}, \cite{N-2}, \cite{RTY}). \\

In the case when $\vert P_{00} + P_{01}\vert = 0$ we can not know the asymptotic profile of the solution $u(t,x)$, and in this case Said-Houari \cite[Theorem 3.3]{said} states that the asymptotic profile still becomes the solution $v(t,x) $ to problem (1.1)-(1.2) with $v_{0} := u_{0} + u_{1}$. He considered such case in terms of the weighted $L^{1}$ data. 

\begin{rem}{\rm If we apply the results due to \cite[Theorem 2.1]{KU} to the "linear" case (i.e., $F(t,x,u) = 0$ in \cite{KU}), their result tells us that
\[t^{\frac{n}{2}(1-\frac{1}{p})}\Vert u(t,\cdot)-(P_{00}+P_{01})G(1+t,\cdot)\Vert_{p} = O(t^{-k/2}),\quad t \to +\infty,\]
for $p \in [1,\infty]$ and $u_{0} \in W_{k}^{1,1} \cap L^{1,\infty}$ and $u_{1} \in L_{k}^{1}\cap L^{\infty}$ with $k \in (0,1]$. So, if we choose $p = 2$, and $k = 1$ in order to compare, we have
\begin{equation}
t^{\frac{n}{4}}\Vert u(t,\cdot)-(P_{00}+P_{01})G(1+t,\cdot)\Vert = O(t^{-1/2}),\quad t \to +\infty.
\end{equation}
The decay order of Theorem 1.2 becomes the same as (1.7). Although we can derive the same type assertion in terms of $L^{\infty}$-norm, too, it is left to the reader's exercise.}
\end{rem}

{\bf Notation.} Throughout this paper, $\| \cdot\|_q$ stands for the usual $L^q({\bf R}^{n})$-norm. For simplicity of notations, in paticular, we use $\| \cdot\|$ instead of $\| \cdot\|_2$. 
Furthermore, we denote the Fourier transform $\hat{\phi}(\xi)$ of the function $\phi(x)$ by
\begin{equation}
{\cal F}(\phi)(\xi) := \hat{\phi}(\xi) := \frac{1}{(2\pi)^{n/2}}\int_{{\bf R}^{n}}e^{-ix\cdot\xi}\phi(x)dx,
\end{equation}
where $i := \sqrt{-1}$, and $x\cdot\xi = \displaystyle{\sum_{i=1}^{n}}x_{i}\xi_{i}$ for $x = (x_{1},\cdots,x_{n})$ and $\xi = (\xi_{1},\cdots,\xi_{n})$, and the inverse Fourier transform of ${\cal F}$ is denoted by ${\cal F}^{-1}$. When we estimate several functions by applying the Fourier transform sometimes we can also use the following definition in place of (1.8)
\[{\cal F}(\phi)(\xi) := \int_{{\bf R}^{n}}e^{-ix\cdot\xi}\phi(x)dx\]
without loss of generality. We also use the notation
\[v_{t}=\frac{\partial u}{\partial t}, \quad v_{tt}=\frac{\partial^{2} v}{\partial t^{2}}, \quad \Delta = \sum^n_{i=1}\frac{\partial^2}{\partial x_i^2},\ \ x=(x_1,\cdots,x_n).\]

\section{Proof of Theorem 1.1.}

In this section, we shall prove Theorem 1.1 by relying on a new method, which has its origin in \cite{Ik-3}.\\

In the following proof we can assume that the initial datum $v_{0}$ are sufficiently smooth, say $v_{0} \in C_{0}^{\infty}({\bf R}_{x}^{n})$ because of the density argument.\\

{\it Proof of Theorem 1.1.}\,First, we apply the Fourier transform of both sides of (1.1)-(1.2), then in the Fourier space ${\bf R}_{\xi}^{n}$ one has the reduced problem:
\begin{equation}
\hat{v}_{t}(t,\xi)+\vert\xi\vert^{2}\hat{v}(t,\xi) = 0,\ \ \ (t,\xi)\in (0,\infty)\times {\bf R}_{\xi}^{n} ,\label{eqn}
\end{equation}
\begin{equation}
\hat{v}(0,\xi)= \hat{v_{0}}(\xi),\,\,\xi \in {\bf R}_{\xi}^{n}.\label{initial}
\end{equation}
Then we can solve (2.1)-(2.2) directly:
\[\hat{v}(t,\xi) = \hat{v_{0}}(\xi) e^{-\vert\xi\vert^{2}t}.\]
We notice that 
\[\hat{v_{0}}(\xi) = \int_{{\bf R}^{n}}v_{0}(x)\cos(x\cdot\xi)dx - i \int_{{\bf R}^{n}}v_{0}(x)\sin(x\cdot\xi)dx\]
\[= \int_{{\bf R}^{n}}v_{0}(x)(\cos(x\cdot\xi)-1)dx - i \int_{{\bf R}^{n}}v_{0}(x)\sin(x\cdot\xi)dx + \int_{{\bf R}^{n}}v_{0}(x)dx\]
\[=: A(\xi) -iB(\xi) + P_{0},\]
so that one has
\[\hat{v}(t,\xi) - P_{0}e^{-\vert\xi\vert^{2}t} = A(\xi)e^{-\vert\xi\vert^{2}t} -iB(\xi)e^{-\vert\xi\vert^{2}t}.\]
This implies
\begin{equation}
\Vert \hat{v}(t,\cdot) - P_{0}e^{-\vert\cdot\vert^{2}t}\Vert^{2} \leq C\int_{{\bf R}^{n}}\vert A(\xi)\vert^{2}e^{-2\vert\xi\vert^{2}t}d\xi + C\int_{{\bf R}^{n}}\vert B(\xi)\vert^{2}e^{-2\vert\xi\vert^{2}t}d\xi.
\end{equation}
Set
\[L := \sup_{\theta\ne 0}\frac{\vert 1-\cos\theta\vert}{\vert\theta\vert} < +\infty,\]
\[M := \sup_{\theta\ne 0}\frac{\vert\sin\theta\vert}{\vert\theta\vert} < +\infty.\]
Then, we can estimate (2.3) as follows: in case of $\xi \ne 0$, for small $\delta > 0$ one has
\[\int_{\vert x\vert \geq \delta}\vert v_{0}(x)\vert\vert(\cos(x\cdot\xi)-1)\vert dx\]
\[\leq \int_{\vert x\vert \geq \delta}\vert v_{0}(x)\vert\frac{\vert(\cos(x\cdot\xi)-1)\vert}{\vert x\cdot\xi\vert}\vert x\cdot\xi\vert dx\]
\begin{equation}
\leq L\vert\xi\vert\int_{{\bf R}^{n}}\vert v_{0}(x)\vert\vert x\vert dx.
\end{equation}
Letting $\delta \downarrow 0$ in (2.4), one has
\begin{equation}
\vert A(\xi)\vert \leq L\vert\xi\vert\int_{{\bf R}^{n}}\vert v_{0}(x)\vert\vert x\vert dx \quad (\xi \in {\bf R}^{n}).
\end{equation} 
Note that (2.5) holds true also in the case when $\xi = 0$. Similarly to (2.5), one also has
\begin{equation}
\vert B(\xi)\vert \leq M\vert\xi\vert\int_{{\bf R}^{n}}\vert v_{0}(x)\vert\vert x\vert dx \quad (\xi \in {\bf R}^{n}).
\end{equation} 
Thus, because of (2.3), (2.5) and (2.6) one can arrive at the meaningful inequality:
\[\Vert \hat{v}(t,\cdot) - P_{0}e^{-\vert\cdot\vert^{2}t}\Vert^{2}\]
\[\leq \C(L^{2}+M^{2})(\int_{{\bf R}^{n}}\vert v_{0}(x)\vert\vert x\vert dx)^{2}\int_{{\bf R}^{n}}\vert\xi\vert^{2}e^{-2\vert\xi\vert^{2}t}d\xi\]
\[\leq C(L^{2}+M^{2})(\int_{{\bf R}^{n}}\vert v_{0}(x)\vert\vert x\vert dx)^{2}t^{-\frac{n}{2}-1}.\]
Finally, because of the Plancherel Theorem and the well-known fact that $G(t,x) = {\cal F}^{-1}(e^{-t\vert\xi\vert^{2}})(x)$, one has the desired estimate.
\hfill
$\Box$

\begin{rem} {\rm 
By observing the proof of Theorem 1.1 we can find that $v_{0} \in L^{1,1}({\bf R}^{n})$ implies $\hat{v}_{0} \in C^{1}({\bf R}^{n})$, and because of the Riemann-Lebesgue theorem one has
\[\lim_{\vert\xi\vert \to +\infty}\frac{\partial\hat{v}_{0}(\xi)}{\partial\xi_{j}} = \lim_{\vert\xi\vert \to +\infty}{\cal F}(-ix_{j}v_{0}(\cdot))(\xi) = 0,\]
\[\vert\frac{\partial\hat{v}_{0}(\xi)}{\partial\xi_{j}}\vert \leq C\Vert v_{0}\Vert_{1,1}\quad (j = 1,2,\cdots,n).\]
Moreover, from the mean value theorem one has
\[\hat{v}_{0}(\xi) - \hat{v}_{0}(0) = \nabla\hat{v}_{0}(\theta\xi)\cdot\xi, \quad \theta \in (0,1),\]
so that one can also arrive at the essential inequality in our proof:
\[\vert \hat{v}_{0}(\xi) - P_{0}\vert \leq C\Vert v_{0}\Vert_{1,1}\vert\xi\vert.\]
Although we can generalize this idea to the initial datum $v_{0} \in L^{1,k}({\bf R}^{n})$ with more heavy weight, this will be our next project.
}
\end{rem}

\section{Proof of Theorem 1.2.}

Let us prove Theorem 1.2 based on an idea due to \cite{Ik-3}. The first part of proof corresponds to the low frequency estimate of the solution. \\
\begin{lem}  It is true that there exists a constant $C>0$ such that for $t > 0$ one has
\[\int_{\vert\xi\vert\leq 1/4}\vert \hat{u}(t,\xi) - (P_{00}+P_{01})e^{-t\vert\xi\vert^{2}}\vert^{2}d\xi\]
\[\leq Ct^{-\frac{n}{2}-1}(\Vert u_{0}\Vert_{1,1}^{2}+\Vert u_{1}\Vert_{1,1}^{2}) + Ct^{-\frac{n}{2}-2}(\Vert u_{0}\Vert_{1}^{2}+\Vert u_{1}\Vert_{1}^{2}) + Ce^{-t}(\Vert u_{0}\Vert^{2} + \Vert u_{1}\Vert^{2}).\]
\end{lem}
{\it Proof of Lemma 3.1.} We apply the Fourier transform of both sides of (1.3)-(1.4), then in the Fourier space ${\bf R}_{\xi}^{n}$ one has the reduced problem:
\begin{equation}
\hat{u}_{tt}(t,\xi)+\vert\xi\vert^{2}\hat{u}(t,\xi) + \hat{u}_{t}(t,\xi) = 0,\ \ \ (t,\xi)\in (0,\infty)\times {\bf R}_{\xi}^{n} ,\label{eqn}
\end{equation}
\begin{equation}
\hat{u}(0,\xi)= \hat{u_{0}}(\xi),\ \ \hat{u}_{t}(0,\xi)= \hat{u_{1}}(\xi),\ \ \ x\in {\bf R}_{\xi}^{n}.\label{initial}
\end{equation}
Let us solve (3.1)-(3.2) directly under the condition that $\vert\xi\vert \leq 1/4$. In this case we get
\begin{equation}
\hat{u}(t,\xi) = \frac{\hat{u_{1}}(\xi)-\sigma_{2}\hat{u_{0}}(\xi)}{\sigma_{1}-\sigma_{2}}e^{\sigma_{1}t} + \frac{\hat{u_{0}}(\xi)\sigma_{1}-\hat{u_{1}}(\xi)}{\sigma_{1}-\sigma_{2}}e^{\sigma_{2}t},
\end{equation}
where $\sigma_{j} \in {\bf R}$ ($j = 1,2$) have a form:
\[\sigma_{1} = \frac{-1+\sqrt{1-4\vert\xi\vert^{2}}}{2}, \quad \sigma_{2} = \frac{-1-\sqrt{1-4\vert\xi\vert^{2}}}{2}.\]
Here, we notice that
\begin{equation}
\sigma_{1}^{2} = -\sigma_{1} - \vert\xi\vert^{2},
\end{equation}
\begin{equation}
\sigma_{2}^{2} = -\sigma_{2} - \vert\xi\vert^{2}.
\end{equation}
By rewriting (3.3) using (3.4) and (3.5) one has
\begin{equation}
\hat{u}(t,\xi) = e^{-t\vert\xi\vert^{2}}\{K_{1}(t,\xi) + K_{2}(t,\xi)\},
\end{equation}
where
\[K_{1}(t,\xi) = \frac{\sigma_{2}\hat{u_{0}}(\xi)-\hat{u_{1}}(\xi)}{\sigma_{2}-\sigma_{1}}e^{-\sigma_{1}^{2}t},\]
\[K_{2}(t,\xi) = \frac{\hat{u_{0}}(\xi)\sigma_{1}-\hat{u_{1}}(\xi)}{\sigma_{1}-\sigma_{2}}e^{-\sigma_{2}^{2}t}.\]
It is important to know that $K_{1}(t,\xi)$ can be decomposed into the following style. This decomposition comes from Chill-Haraux \cite{CH}. 
\[K_{1}(t,\xi) = \hat{u_{0}}(\xi)+\hat{u_{1}}(\xi) + \frac{\sigma_{1}\hat{u_{0}}(\xi)}{\sigma_{2}-\sigma_{1}}\] 
\[+ \frac{\sigma_{2}\hat{u_{0}}(\xi)(1-e^{-\sigma_{1}^{2}t})}{\sigma_{1}-\sigma_{2}} + \frac{\hat{u_{1}}(\xi)(e^{-\sigma_{1}^{2}t}-(\sigma_{1}-\sigma_{2}))}{\sigma_{1}-\sigma_{2}}.\]
So one has arrived at the meaningful relation:
\[\hat{u}(t,\xi) = e^{-t\vert\xi\vert^{2}}\{\hat{u_{0}}(\xi)+\hat{u_{1}}(\xi) + \frac{\sigma_{1}\hat{u_{0}}(\xi)}{\sigma_{2}-\sigma_{1}}\]
\begin{equation}
+  \frac{\sigma_{2}\hat{u_{0}}(\xi)(1-e^{-\sigma_{1}^{2}t})}{\sigma_{1}-\sigma_{2}} + \frac{\hat{u_{1}}(\xi)(e^{-\sigma_{1}^{2}t}-(\sigma_{1}-\sigma_{2}))}{\sigma_{1}-\sigma_{2}}\} + e^{-t\vert\xi\vert^{2}}K_{2}(t,\xi).
\end{equation}
Now let us use the idea similar to the proof of Theorem 1.1 (see \cite[Lemma 3.1]{Ik-3}). We use the relations
\begin{equation}
\hat{u}_{j}(\xi) = A_{j}(\xi) -iB_{j}(\xi) + P_{0j}\quad (j = 0,1), 
\end{equation}
where
\[A_{j}(\xi) := \int_{{\bf R}^{n}}(\cos(x\cdot\xi)-1)u_{j}(x) dx, \quad B_{j}(\xi) := \int_{{\bf R}^{n}}\sin(x\cdot\xi)u_{j}(x) dx, \quad(j = 0,1).\]
Because of (3.7) and (3.8) we get the useful identity for all $\xi$ satisfying $\vert\xi\vert \leq 1/4$:
\[\hat{u}(t,\xi)\ - (P_{00}+P_{01})e^{-t\vert\xi\vert^{2}}\]
\begin{equation}
= (A_{0}(\xi)-iB_{0}(\xi) + A_{1}(\xi)-iB_{1}(\xi))e^{-t\vert\xi\vert^{2}}
\end{equation}
\begin{equation}
+ e^{-t\vert\xi\vert^{2}}\{\frac{\sigma_{1}\hat{u_{0}}(\xi)}{\sigma_{2}-\sigma_{1}} +  \frac{\sigma_{2}\hat{u_{0}}(\xi)(1-e^{-\sigma_{1}^{2}t})}{\sigma_{1}-\sigma_{2}}+ \frac{\hat{u_{1}}(\xi)(e^{-\sigma_{1}^{2}t}-(\sigma_{1}-\sigma_{2}))}{\sigma_{1}-\sigma_{2}}\} + e^{-t\vert\xi\vert^{2}}K_{2}(t,\xi).
\end{equation}
The essential part of our result is in the estimation for (3.9). Although the estimates for (3.10) can be done almost similar to \cite{CH}, for the sake of completeness of the proof we will write down all estimates for (3.9) and (3.10).\\

{\bf (I)\,The $L^{2}$ estimates for (3.10).}\\ 

First, let us note the inequalities that
\begin{equation}
0 < -\sigma_{1} \leq 2\vert\xi\vert^{2},
\end{equation} 
\begin{equation}
\sigma_{1}-\sigma_{2} = \sqrt{1-4\vert\xi\vert^{2}} \geq \sqrt{3/4} > 1/2.
\end{equation} 
So, because of (3.11) and (3.12), the first term of (3.10) can be estimated in terms of $L^{2}$-norm as follows:
\[J_{0}(t) := \int_{\vert\xi\vert \leq 1/4}e^{-2t\vert\xi\vert^{2}}\vert\frac{\sigma_{1}\hat{u_{0}}(\xi)}{\sigma_{2}-\sigma_{1}}\vert^{2}d\xi\]
\begin{equation}
\leq C\Vert \hat{u}_{0}\Vert_{\infty}^{2}\int_{\vert\xi\vert \leq 1/4}e^{-2t\vert\xi\vert^{2}}\vert\xi\vert^{4}d\xi \leq Ct^{-\frac{n}{2}-2}\Vert u_{0}\Vert_{1}^{2}.
\end{equation}
Furthermore, because of the mean value theorem we know
\[\vert 1-e^{-\sigma_{1}^{2}t}\vert \leq t\sigma_{1}^{2},\]
so that one has
\[J_{1}(t):= \int_{\vert\xi\vert\leq 1/4}\frac{\sigma_{2}^{2}(1-e^{-\sigma_{1}^{2}t})^{2}}{\vert \sigma_{1}-\sigma_{2}\vert^{2}}\vert\hat{u}_{0}(\xi)\vert^{2}e^{-2t\vert\xi\vert^{2}}d\xi\]
\[\leq 4t^{2}\Vert u_{0}\Vert_{1}^{2}\int_{\vert\xi\vert\leq 1/4}\sigma_{2}^{2}\sigma_{1}^{4}e^{-2t\vert\xi\vert^{2}}d\xi.\]
Since $\sigma_{1}^{4} \leq 16\vert\xi\vert^{8}$, and
\[\sigma_{2}^{2} = -\sigma_{2}-\vert\xi\vert^{2} \leq -\sigma_{2} = \frac{1+\sqrt{1-4\vert\xi\vert^{2}}}{2} \leq \frac{1}{2} + \frac{1}{2} = 1,\]
one gets
\[\int_{\vert\xi\vert\leq 1/4}\sigma_{2}^{2}\sigma_{1}^{4}e^{-2t\vert\xi\vert^{2}}d\xi \leq 16\int_{\vert\xi\vert\leq 1/4}\vert\xi\vert^{8}e^{-2t\vert\xi\vert^{2}}d\xi \leq Ct^{-\frac{n}{2}-4},\]
so one has
\begin{equation}
J_{1}(t) \leq Ct^{-\frac{n}{2}-2}\Vert u_{0}\Vert_{1}^{2}.
\end{equation}
On the other hand, since
\[\vert e^{-\sigma_{1}^{2}t}-(\sigma_{1}-\sigma_{2})\vert = \vert(e^{-\sigma_{1}^{2}t}-1) + 1- (\sigma_{1}-\sigma_{2}) \vert\]
\[= \vert (e^{-\sigma_{1}^{2}t}-1)+ (1-\sqrt{1-4\vert\xi\vert^{2}})\vert =  \vert(e^{-\sigma_{1}^{2}t}-1) + \frac{4\vert\xi\vert^{2}}{(1+\sqrt{1-4\vert\xi\vert^{2}})}\vert \]
\[\leq \vert e^{-\sigma_{1}^{2}t}-1\vert  + 4\vert\xi\vert^{2},\]
because of the mean value theorem again, we see that
\begin{equation}
\vert e^{-\sigma_{1}^{2}t}-(\sigma_{1}-\sigma_{2})\vert \leq t\sigma_{1}^{2} + 4\vert\xi\vert^{2} \leq 4t\vert\xi\vert^{4} + 4\vert\xi\vert^{2}.
\end{equation}
Thus, from (3.15) one can estimate as follows:
\[J_{2}(t) := \int_{\vert\xi\vert\leq1/4}e^{-2t\vert\xi\vert^{2}}\frac{\vert\hat{u_{1}}(\xi)\vert^{2}\vert e^{-\sigma_{1}^{2}t}-(\sigma_{1}-\sigma_{2}))\vert^{2}}{\vert\sigma_{1}-\sigma_{2}\vert^{2}}d\xi\]
\begin{equation}
\leq C \Vert u_{1}\Vert_{1}^{2}\int_{\vert\xi\vert\leq1/4}e^{-2t\vert\xi\vert^{2}}(t^{2}\vert\xi\vert^{8} + \vert\xi\vert^{4})d\xi \leq C \Vert u_{1}\Vert_{1}^{2}t^{-\frac{n}{2}-2}.
\end{equation}

{\bf (II)\,The $L^{2}$ estimates for (3.9).}\\ 
Let us estimate (3.9) in terms of $L^{2}$-norm, which is the main part of our result.\\ 
In fact,
\[J_{3}(t) := \int_{\vert\xi\vert\leq 1/4} \vert A_{0}(\xi)-iB_{0}(\xi) + A_{1}(\xi)-iB_{1}(\xi)\vert^{2}e^{-2t\vert\xi\vert^{2}}d\xi\]
\begin{equation}
\leq C\int_{\vert\xi\vert\leq 1/4}(\vert A_{0}(\xi)\vert + \vert B_{0}(\xi)\vert + \vert A_{1}(\xi)\vert + \vert B_{1}(\xi)\vert)^{2}e^{-2t\vert\xi\vert^{2}}d\xi.
\end{equation}
By proceeding the same computations as in the proof of Theorem 1.1 one can get
\[\vert A_{j}(\xi)\vert \leq L\vert\xi\vert\Vert u_{j}\Vert_{1,1},\quad (j = 0,1)\]
\[\vert B_{j}(\xi)\vert \leq M\vert\xi\vert\Vert u_{j}\Vert_{1,1},\quad (j = 0,1),\]
where we have set again
\[L := \sup_{\theta\ne 0}\frac{\vert 1-\cos\theta\vert}{\vert\theta\vert} < +\infty, \quad M := \sup_{\theta\ne 0}\frac{\vert\sin\theta\vert}{\vert\theta\vert} < +\infty,\]
so that we have
\[J_{3}(t) \leq C(\Vert u_{0}\Vert_{1,1}^{2}+\Vert u_{1}\Vert_{1,1}^{2})\int_{\vert\xi\vert\leq 1/4}(L+M)^{2}\vert\xi\vert^{2}e^{-2t\vert\xi\vert^{2}}d\xi \]
\begin{equation}
\leq C(L+M)^{2}(\Vert u_{0}\Vert_{1,1}^{2}+\Vert u_{1}\Vert_{1,1}^{2})t^{-\frac{n}{2}-1}.
\end{equation}

{\bf (III)\,The $L^{2}$ estimates for the last term of (3.10).}\\ 
Let us estimate the following $J_{4}(t)$: 
\[J_{4}(t) := \int_{\vert\xi\vert \leq 1/4}e^{-2t\vert\xi\vert^{2}}\vert K_{2}(t,\xi)\vert^{2}d\xi\]
\[= \int_{\vert\xi\vert \leq 1/4}e^{-2t\vert\xi\vert^{2}}\vert\frac{\hat{u}_{0}(\xi)\sigma_{1}-\hat{u}_{1}(\xi)}{\sigma_{1}-\sigma_{2}} \vert^{2}e^{-2\sigma_{2}^{2}t}d\xi.\]
Since $\sigma_{1}^{2} \leq -\sigma_{1} \leq 1/8$ because of (3.11) and
\[\sigma_{2}^{2} = \frac{1+\sqrt{1-4\vert\xi\vert^{2}}}{2} - \vert\xi\vert^{2} \geq \frac{1+\sqrt{1-4/16}}{2}-\frac{1}{16} = \frac{7+4\sqrt{3}}{16} > \frac{1}{2},\]
by (3.12) and the Plancherel theorem one has
\[J_{4}(t) \leq Ce^{-t}\int_{\vert\xi\vert \leq 1/4}e^{-2t\vert\xi\vert^{2}}(\vert\hat{u}_{0}(\xi)\vert^{2}\sigma_{1}^{2} + \vert\hat{u}_{1}(\xi)\vert^{2})d\xi\]
\begin{equation}
\leq Ce^{-t}\int_{\vert\xi\vert \leq 1/4}(\vert\hat{u}_{0}(\xi)\vert^{2} + \vert\hat{u}_{1}(\xi)\vert^{2})d\xi \leq Ce^{-t}(\Vert u_{0}\Vert^{2} + \Vert u_{1}\Vert^{2}).
\end{equation}

Finally, because of (3.9), (3.10), (3.13), (3.14), (3.16), (3.18) and (3.19) one has arrived at the desired estimate for lemma 3.1:
\[\int_{\vert\xi\vert\leq 1/4}\vert \hat{u}(t,\xi) - (P_{00}+P_{01})e^{-t\vert\xi\vert^{2}}\vert^{2}d\xi\]
\[\leq Ct^{-\frac{n}{2}-1}(\Vert u_{0}\Vert_{1,1}^{2}+\Vert u_{1}\Vert_{1,1}^{2}) + Ct^{-\frac{n}{2}-2}(\Vert u_{0}\Vert_{1}^{2}+\Vert u_{1}\Vert_{1}^{2}) + Ce^{-t}(\Vert u_{0}\Vert^{2} + \Vert u_{1}\Vert^{2}).\]
\hfill
$\Box$

\par

The second part of proof corresponds to the high frequency estimate of the solution.
\begin{lem}  It is true that there exists a constant $C>0$ such that for $t > 0$ one has
\[\int_{\vert\xi\vert\geq 1}\vert \hat{u}(t,\xi) - (P_{00}+P_{01})e^{-t\vert\xi\vert^{2}}\vert^{2}d\xi\]
\[\leq Ct^{-\frac{n}{2}-1}(\Vert u_{0}\Vert_{1,1}^{2}+\Vert u_{1}\Vert_{1,1}^{2}) + Ce^{-t}(\Vert\nabla u_{0}\Vert^{2}+\Vert u_{1}\Vert^{2}).\]
\end{lem}
{\it Proof of Lemma 3.2.} We start with the following explicit formula under the assumption $\vert\xi\vert \geq 1$:
\begin{equation}
\hat{u}(t,\xi) = \frac{\hat{u_{1}}(\xi)-\lambda_{2}\hat{u_{0}}(\xi)}{\lambda_{1}-\lambda_{2}}e^{\lambda_{1}t} + \frac{\hat{u_{0}}(\xi)\lambda_{1}-\hat{u_{1}}(\xi)}{\lambda_{1}-\lambda_{2}}e^{\lambda_{2}t},
\end{equation}
where $\lambda_{j} \in {\bf C}$ ($j = 1,2$) have a form:
\[\lambda_{1} = \frac{-1+\sqrt{4\vert\xi\vert^{2}-1}i}{2}, \quad \lambda_{2} = \frac{-1-\sqrt{4\vert\xi\vert^{2}-1}i}{2}.\]
Here, we also notice that
\begin{equation}
\lambda_{j}^{2} = -\lambda_{j} - \vert\xi\vert^{2},\quad (j = 1,2).
\end{equation}
By rewriting (3.20) using (3.21) one has
\begin{equation}
\hat{u}(t,\xi) = e^{-t\vert\xi\vert^{2}}\{\frac{\lambda_{2}\hat{u_{0}}(\xi)-\hat{u_{1}}(\xi)}{\lambda_{2}-\lambda_{1}}e^{-\lambda_{1}^{2}t} + \frac{\hat{u_{0}}(\xi)\lambda_{1}-\hat{u_{1}}(\xi)}{\lambda_{1}-\lambda_{2}}e^{-\lambda_{2}^{2}t}\}.
\end{equation}
We set for later use. 
\[K_{3}(t,\xi) = \frac{\hat{u_{0}}(\xi)\lambda_{1}-\hat{u_{1}}(\xi)}{\lambda_{1}-\lambda_{2}}e^{-\lambda_{2}^{2}t}.\]
By the same procedure as in (3.7), (3.8), (3.9) and (3.10) one has the following decomposition again: 
for all $\xi$ satisfying $\vert\xi\vert \geq 1$:
\[\hat{u}(t,\xi)\ - (P_{00}+P_{01})e^{-t\vert\xi\vert^{2}}\]
\begin{equation}
= (A_{0}(\xi)-iB_{0}(\xi) + A_{1}(\xi)-iB_{1}(\xi))e^{-t\vert\xi\vert^{2}}
\end{equation}
\begin{equation}
+ e^{-t\vert\xi\vert^{2}}\{\frac{\lambda_{1}\hat{u_{0}}(\xi)}{\lambda_{2}-\lambda_{1}} +  \frac{\lambda_{2}\hat{u_{0}}(\xi)(1-e^{-\lambda_{1}^{2}t})}{\lambda_{1}-\lambda_{2}}+ \frac{\hat{u_{1}}(\xi)(e^{-\lambda_{1}^{2}t}-(\lambda_{1}-\lambda_{2}))}{\lambda_{1}-\lambda_{2}}\} + e^{-t\vert\xi\vert^{2}}K_{3}(t,\xi).
\end{equation}
\noindent
In the computations below we have to use the following relation for $\lambda_{j} \in {\bf C}$ ($j = 1,2$)
\begin{equation}
\vert\lambda_{j}\vert = \vert\xi\vert \geq 1, \quad (j = 1,2)
\end{equation} 
\begin{equation}
\vert\lambda_{1}-\lambda_{2}\vert = \sqrt{4\vert\xi\vert^{2}-1} \geq \sqrt{3}.
\end{equation} 
\begin{equation}
\lambda_{1}^{2} = \frac{1}{2}-\vert\xi\vert^{2} - \frac{\sqrt{4\vert\xi\vert^{2}-1}i}{2},
\end{equation}
\begin{equation}
\lambda_{2}^{2} = \frac{1}{2}-\vert\xi\vert^{2} + \frac{\sqrt{4\vert\xi\vert^{2}-1}i}{2},
\end{equation}
so that one has
\begin{equation}
\vert e^{-\lambda_{j}^{2}t}\vert = e^{-t/2}e^{\vert\xi\vert^{2}t},\quad (j = 1,2).
\end{equation}

{\bf (I)\,The $L^{2}$ estimates for (3.24).}\\ 

The first term of (3.24) can be estimated in terms of $L^{2}$-norm as follows, because of (3.25) and (3.26). In this case we need the regularity on the initial amplitude $u_{0}$.
\[I_{0}(t) := \int_{\vert\xi\vert \geq 1}e^{-2t\vert\xi\vert^{2}}\vert\frac{\lambda_{1}\hat{u_{0}}(\xi)}{\lambda_{2}-\lambda_{1}}\vert^{2}d\xi\]
\begin{equation}
\leq \frac{e^{-2t}}{3}\int_{\vert\xi\vert \geq 1}\vert\xi\vert^{2}\vert\hat{u}_{0}(\xi)\vert^{2}d\xi \leq Ce^{-2t}\Vert\nabla u_{0}\Vert^{2}.
\end{equation}
Furthermore, it follows from (3.25), (3.26) and (3.29) that
\[I_{1}(t):= \int_{\vert\xi\vert\geq 1}\frac{\vert\lambda_{2}\vert^{2}\vert1-e^{-\lambda_{1}^{2}t}\vert^{2}}{\vert \lambda_{1}-\lambda_{2}\vert^{2}}\vert\hat{u}_{0}(\xi)\vert^{2}e^{-2t\vert\xi\vert^{2}}d\xi\]
\[\leq C\int_{\vert\xi\vert\geq 1}\frac{\vert\lambda_{2}\vert^{2}}{\vert \lambda_{1}-\lambda_{2}\vert^{2}}\vert\hat{u}_{0}(\xi)\vert^{2}e^{-2t\vert\xi\vert^{2}}d\xi +C\int_{\vert\xi\vert\geq 1}\frac{\vert\lambda_{2}\vert^{2}}{\vert \lambda_{1}-\lambda_{2}\vert^{2}}\vert e^{-\lambda_{1}^{2}t}\vert^{2}\vert\hat{u}_{0}(\xi)\vert^{2}e^{-2t\vert\xi\vert^{2}}d\xi \]
\[\leq \frac{e^{-2t}}{3}\int_{\vert\xi\vert\geq1}\vert\xi\vert^{2}\vert\hat{u}_{0}(\xi)\vert^{2}d\xi + \frac{e^{-t}}{3}\int_{\vert\xi\vert\geq1}e^{-2t\vert\xi\vert^{2}}\vert\xi\vert^{2}\vert\hat{u}_{0}(\xi)\vert^{2}e^{2t\vert\xi\vert^{2}}d\xi\]
\begin{equation}
\leq Ce^{-2t}\Vert\nabla u_{0}\Vert^{2} + Ce^{-t}\Vert\nabla u_{0}\Vert^{2}.
\end{equation}
\noindent
On the other hand, because of (3.26) and (3.29) one gets
\[I_{2}(t) := \int_{\vert\xi\vert\geq1}e^{-2t\vert\xi\vert^{2}}\vert\frac{\hat{u_{1}}(\xi)(e^{-\lambda_{1}^{2}t}-(\lambda_{1}-\lambda_{2}))}{\lambda_{1}-\lambda_{2}}\vert ^{2}d\xi\]
\[\leq C \int_{\vert\xi\vert \geq 1}\frac{\vert e^{-\lambda_{1}^{2}t}\vert^{2}}{\vert\lambda_{1}-\lambda_{2}\vert^{2}}e^{-2t\vert\xi\vert^{2}}\vert\hat{u}_{1}(\xi)\vert^{2}d\xi 
+ C \int_{\vert\xi\vert \geq 1}e^{-2t\vert\xi\vert^{2}}\vert\hat{u}_{1}(\xi)\vert^{2}d\xi\]
\begin{equation}
\leq C \int_{\vert\xi\vert \geq 1}e^{-2t\vert\xi\vert^{2}}e^{-t}e^{2t\vert\xi\vert^{2}}\vert\hat{u}_{1}(\xi)\vert^{2}d\xi + C e^{-2t}\int_{\vert\xi\vert \geq 1}\vert\hat{u}_{1}(\xi)\vert^{2}d\xi \leq Ce^{-t}\Vert u_{1}\Vert^{2}.
\end{equation}

{\bf (II)\,The $L^{2}$ estimates for the last term of (3.24).}\\ 
Let us estimate the following $I_{3}(t)$ by using (3.25), (3.26) and (3.29): 
\[I_{3}(t) := \int_{\vert\xi\vert \geq 1}e^{-2t\vert\xi\vert^{2}}\vert K_{3}(t,\xi)\vert^{2}d\xi = \int_{\vert\xi\vert \geq 1}e^{-2t\vert\xi\vert^{2}}\vert\frac{\hat{u}_{0}(\xi)\lambda_{1}-\hat{u}_{1}(\xi)}{\lambda_{1}-\lambda_{2}} \vert^{2}\vert e^{-\lambda_{2}^{2}t}\vert^{2}d\xi.\]
\[\leq \frac{1}{3}\int_{\vert\xi\vert \geq 1}e^{-2t\vert\xi\vert^{2}}\vert\lambda_{1}\vert^{2}\vert\hat{u}_{0}(\xi)\vert^{2}\vert e^{-\lambda_{2}^{2}t}\vert^{2}d\xi + \frac{1}{3}\int_{\vert\xi\vert \geq 1}e^{-2t\vert\xi\vert^{2}}\vert\hat{u}_{1}(\xi)\vert^{2}\vert e^{-\lambda_{2}^{2}t}\vert^{2}d\xi\]
\[\leq C\int_{\vert\xi\vert \geq 1}\vert\xi\vert^{2}e^{-t}e^{2t\vert\xi\vert^{2}}e^{-2t\vert\xi\vert^{2}}\vert\hat{u}_{0}(\xi)\vert^{2}d\xi + C\int_{\vert\xi\vert \geq 1}e^{-t}e^{2t\vert\xi\vert^{2}}e^{-2t\vert\xi\vert^{2}}\vert\hat{u}_{1}(\xi)\vert^{2}d\xi\]
\begin{equation}
\leq Ce^{-t}\Vert\nabla u_{0}\Vert^{2} + Ce^{-t}\Vert u_{1}\Vert^{2}.
\end{equation}

{\bf (III)\,The $L^{2}$ estimates for (3.23).}\\ 
Let us estimate (3.23) in terms of $L^{2}$-norm. This part is treated with the similar procedure as in that of (3.18).\\ 
In fact,
\[I_{4}(t) := \int_{\vert\xi\vert\geq 1} \vert A_{0}(\xi)-iB_{0}(\xi) + A_{1}(\xi)-iB_{1}(\xi)\vert^{2}e^{-2t\vert\xi\vert^{2}}d\xi\]
\begin{equation}
\leq C(L\Vert u_{0}\Vert_{1,1}+M\Vert u_{1}\Vert_{1,1} )^{2}\int_{\vert\xi\vert\geq 1}\vert\xi\vert^{2}e^{-2t\vert\xi\vert^{2}}d\xi,
\end{equation}
where we have set again
\[L := \sup_{\theta\ne 0}\frac{\vert 1-\cos\theta\vert}{\vert\theta\vert} < +\infty, \quad M := \sup_{\theta\ne 0}\frac{\vert\sin\theta\vert}{\vert\theta\vert} < +\infty.\]
so that we have
\[I_{4}(t) \leq C(\Vert u_{0}\Vert_{1,1}^{2}+\Vert u_{1}\Vert_{1,1}^{2})\int_{\vert\xi\vert\geq 1}\vert\xi\vert^{2}e^{-2t\vert\xi\vert^{2}}d\xi \]
\begin{equation}
\leq C(\Vert u_{0}\Vert_{1,1}^{2}+\Vert u_{1}\Vert_{1,1}^{2})t^{-\frac{n}{2}-1} \quad (t > 0).
\end{equation}
\noindent
Finally, because of (3.23), (3.24), (3.30), (3.31), (3.32), (3.33) and (3.35) one has arrived at the desired estimate:
\[\int_{\vert\xi\vert\geq 1}\vert \hat{u}(t,\xi) - (P_{00}+P_{01})e^{-t\vert\xi\vert^{2}}\vert^{2}d\xi\]
\[\leq Ct^{-\frac{n}{2}-1}(\Vert u_{0}\Vert_{1,1}^{2}+\Vert u_{1}\Vert_{1,1}^{2}) + Ce^{-t}(\Vert\nabla u_{0}\Vert^{2}+\Vert u_{1}\Vert).\]
\hfill
$\Box$

\par    
{\it Proof of Theorem 1.2.}\,Under the preparation from Lemmas 3.1 and 3.2, we can prove Theorem 1.2 as follows.

In fact, we first make a decomposition as follows by relying on the Plancherel theorem:  
\[\Vert u(t,\cdot)-(P_{00}+P_{01})G(t,\cdot))\Vert^{2},\]
\[\leq (\int_{\vert\xi\vert \leq 1/4} + \int_{1/4 \leq \vert\xi\vert \leq 1} + \int_{1 \leq \vert\xi\vert})\vert \hat{u}(t,\xi)-(P_{00}+P_{01})e^{-t\vert\xi\vert^{2}}\vert^{2}d\xi\]
\[=: R_{1}(t) + R_{2}(t) + R_{3}(t).\]
We can rely on Lemmas 3.1 and 3.2 to get
\[R_{1}(t) \leq Ct^{-\frac{n}{2}-1}(\Vert u_{0}\Vert_{1,1}^{2}+\Vert u_{1}\Vert_{1,1}^{2}) + Ct^{-\frac{n}{2}-2}(\Vert u_{0}\Vert_{1}^{2}+\Vert u_{1}\Vert_{1}^{2}) + Ce^{-t}(\Vert u_{0}\Vert^{2} + \Vert u_{1}\Vert^{2}),\] 
\[R_{3}(t) \leq Ct^{-\frac{n}{2}-1}(\Vert u_{0}\Vert_{1,1}^{2}+\Vert u_{1}\Vert_{1,1}^{2}) + Ce^{-t}(\Vert\nabla u_{0}\Vert^{2}+\Vert u_{1}\Vert^{2}).\]
On the other hand, concerning $R_{2}(t)$ we can apply the classical estimates prepared by Matsumura \cite{Ma} (see also \cite{Ik-3}).  
\[R_{2}(t) = \int_{1/4 \leq \vert\xi\vert \leq 1}\vert \hat{u}(t,\xi)-(P_{00}+P_{01})e^{-t\vert\xi\vert^{2}}\vert^{2}d\xi\]
\[\leq C\int_{1/4 \leq \vert\xi\vert \leq 1}\vert \hat{u}(t,\xi)\vert^{2}d\xi + \vert P_{00}+P_{01}\vert^{2}\int_{1/4 \leq \vert\xi\vert \leq 1}e^{-2t\vert\xi\vert^{2}}d\xi\]
\[\leq Cte^{-t}(\Vert u_{1}\Vert^{2}+\Vert u_{0}\Vert^{2}) + Cte^{-(1-\sqrt{3}/2)t}(\Vert u_{1}\Vert^{2} + \Vert u_{0}\Vert^{2})\]
\[+ C_{n}(\Vert u_{0}\Vert_{1}^{2}+\Vert u_{1}\Vert_{1}^{2})e^{-t/8}.\]
with $C_{n} := \displaystyle{\int_{1/4\leq\vert\xi\vert \leq 1}}d\xi$. These estimates imply the desired statement of Theorem 1.2.
\hfill
$\Box$

\par

\noindent{\em Acknowledgement.}
\smallskip
The work of the author (R. IKEHATA) was supported in part by Grant-in-Aid for Scientific Research (C)22540193 and (A)22244009 of JSPS.


\end{document}